\documentclass[12pt]{article}
\usepackage{graphicx}
\usepackage{amsmath,amsfonts}
\usepackage{mathrsfs}

\let\scr\mathscr
\let\goth\mathfrak

\def\zs#1{_{\lower 3pt \hbox{$\scriptstyle#1$}}}
\def\Pb{\mathbf{P}}
\def\Ex{\mathbf{E}}

\def\onepicsize{13cm}
\def\onepic#1#2{\nobreak\bigskip\noindent\centerline{%
\hbox to \onepicsize{\vbox{%
\hbox to \onepicsize{\includegraphics*[width=\onepicsize]{#1}}%
\hbox to \onepicsize{\hfill\small #2\hfill}}}}\bigskip}

\begin{document}
\date{}
\title{On the Goodness-of-Fit Tests  for  Some Continuous Time Processes}
\author{Sergue\"{\i} Dachian and Yury A. Kutoyants\\
\medskip
{\it Laboratoire de Math\'ematiques,  Universit\'e Blaise
Pascal}\\
{\it Laboratoire de Statistique et Processus, Universit\'e du
Maine}}

\maketitle

%\markboth{S. Dachian and Y. A. Kutoyants}{Sampling in survival analysis}

 \begin{abstract} We present a review of several results
concerning the construction of the Cram\'er-von Mises and
Kolmogorov-Smirnov type goodness-of-fit tests for continuous time
processes. As the models we take a stochastic differential
equation with small noise, ergodic diffusion process, Poisson
process and self-exciting point processes. For every model we
propose the tests which provide the asymptotic size $\alpha $ and
discuss the behaviour of the power function under local
alternatives.  The results of numerical simulations of the tests
are presented.
 \end{abstract}

\bigskip


\vspace*{30pt}

\hrule

\section{Introduction}

The goodness-of-fit tests play an important role in the classical mathematical
statistics. Particularly, the tests of Cram\'er-von Mises, Kolmogorov-Smirnov
and Chi-Squared are well studied and allow to verify the correspondence of the
mathematical models to the observed data (see, for example, Durbin (1973) or
Greenwood and Nikulin (1996)). The similar problem, of course, exists for the
continuous time stochastic processes. The diffusion and Poisson processes are
widely used as mathematical models of many evolution processes in Biology,
Medicine, Physics, Financial Mathematics and in many others fields.  For
example, some theory can propose a diffusion process
$$
{\rm d}X_t=S_*\left(X_t\right)\,{\rm d}t+\sigma \;{\rm d}W_t,\quad X_0, \quad
0\leq t\leq T
$$
as an appropriate model for description of the real data $\left\{X_t,0\leq
t\leq T\right\}$ and we can try to construct an algorithm to verify if this
model corresponds well to these data. The model here is totally defined by the
trend coefficient $S_*\left(\cdot \right)$, which is supposed (if the theory
is true) to be known.  We do not discuss here the problem of verification if
the process $\left\{W_t,0\leq t\leq T\right\}$ is Wiener. This problem is much
more complicated and we suppose that the {\sl noise is white Gaussian}.
Therefore we have a basic hypothesis defined by the trend coefficient
$S_*\left(\cdot \right)$ and we have to test this hypothesis against any other
alternative. Any other means that the observations come from stochastic
differential equation
$$
{\rm d}X_t=S\left(X_t\right)\,{\rm d}t+\sigma \;{\rm d}W_t,\quad X_0, \quad
0\leq t\leq T ,
$$
where $S\left(\cdot \right)\neq S_*\left(\cdot \right)$. We propose some tests
which are in some sense similar to the Cram\'er-von Mises and Kolmogorov-Smirnov
tests. The advantage of classical tests is that they are distribution-free,
i.e., the distribution of the underlying statistics do not depend on the basic
model and this property allows to choose the {\sl universal thresholds}, which
can be used for all models.

For example, if we observe $n$ independent identically distributed random
variables $\left(X_1,\ldots,X_n\right)=X^n$ with distribution function
$F\left(x\right)$ and the basic hypothesis is simple : $F\left(x\right)\equiv
F_*\left(x\right)$, then the Cram\'er-von Mises $W_n^2$ and Kolmogorov-Smirnov
$D_n$ statistics are
$$
W_n^2=n\int_{-\infty }^{\infty }\left[\hat
F_n\left(x\right)-F_*\left(x\right)\right]^2 \,{\rm d}F_*\left(x\right),\qquad
\quad D_n=\sup_x\left|\hat F_n\left(x\right)-F_*\left(x\right) \right|
$$
respectively. Here
$$
\hat F_n\left(x\right)=\frac{1}{n}\sum_{j=1}^{n}1_{\left\{X_j<x\right\}}
$$
is the empirical distribution function.  Let us denote by
$\left\{W_0\left(s\right), 0\leq s\leq 1\right\}$  a Brownian bridge, i.e., a
continuous Gaussian
process with
$$
\Ex W_0\left(s\right)=0,\qquad \Ex W_0\left(s\right)W_0\left(t\right)= t\wedge
s-st.
$$
Then the
limit behaviour of these statistics  can be described with the help of this
process as follows
$$
W_n^2\Longrightarrow \int_{0}^{1}W_0\left(s\right)^2{\rm d}s,\qquad \quad
\sqrt{n}D_n \Longrightarrow \sup_{0\leq s\leq
1}\left|W_0\left(s\right)\right|.
$$
Hence the corresponding Cram\'er-von Mises and Kolmogorov-Smirnov tests
$$
\psi_n\left(X^n\right)=1_{\left\{W _n^2>c_\alpha \right\}},\qquad
\phi_n\left(X^n\right)=1_{\left\{\sqrt{n}D_n>d_\alpha \right\}}
$$
with constants $c_\alpha , d_\alpha  $ defined by the equations
$$
\Pb\left\{\int_{0}^{1}W_0\left(s\right)^2{\rm d}s > c_\alpha \right\}=\alpha
,\qquad \Pb\left\{\sup_{0\leq s\leq 1}\left|W_0\left(s\right)\right|> d_\alpha
\right\}=\alpha
$$
are of asymptotic size $\alpha $.  It is easy to see that these tests are
distribution-free (the limit distributions do not depend of the function
$F_*\left(\cdot \right)$) and are consistent against any fixed alternative
(see, for example, Durbin (1973)).

It is interesting to study these tests for {\sl nondegenerate set of
alternatives}, i.e., for alternatives with limit power function less than
1. It can be realized on the close nonparametric alternatives of the special
form making this problem asymptotically equivalent to the {\sl signal in
Gaussian noise} problem. Let us put
$$
F\left(x\right)=F_*\left(x\right)+\frac{1}{\sqrt{n}}\int_{-\infty
}^{x}h\left(F_*\left(y\right)\right)\;{\rm d}F_*\left(y\right) ,
$$
where the function $h\left(\cdot \right)$ describes the alternatives.
We suppose that
$$
\int_{0}^{1}h\left(s\right)\;{\rm d}s=0,\qquad
\int_{0}^{1}h\left(s\right)^2\;{\rm d}s<\infty .
$$
Then we
have the following convergence (under fixed alternative, given by the function
$h\left(\cdot \right)$):
\begin{align*}
&W_n^2\Longrightarrow \int_{0}^{1}\left[\int_{0}^{s}h\left(v\right){\rm
d}v+W_0\left(s\right)\right]^2{\rm d}s ,\\
& \sqrt{n}D_n \Longrightarrow
\sup_{0\leq s\leq 1}\left|\int_{0}^{s}h\left(v\right){\rm
d}v+W_0\left(s\right)\right|
\end{align*}
We see that this problem is asymptotically equivalent to the  following {\sl
signal in Gaussian noise} problem:
\begin{equation}
\label{sn}
{\rm d}Y_s=h_*\left(s\right)\,{\rm d}s+{\rm d}W_0\left(s\right),\quad 0\leq
s\leq 1.
\end{equation}

Indeed, if we use the statistics
$$
W^2=\int_{0}^{1}Y_s^2\;{\rm d}s,\qquad D=\sup_{0\leq s\leq 1}\left|Y_s\right|
$$
then under hypothesis $h\left(\cdot \right)\equiv 0$ and alternative
$h\left(\cdot \right)\neq  0$   the distributions of these
statistics coincide with the limit distributions of $W_n^2$ and
$\sqrt{n}D_n$ under hypothesis and alternative respectively.

Our goal is to see how such kind of tests can be constructed in the case of
continuous time models of observation and particularly in the cases of some
diffusion and point processes.  We consider the diffusion processes with small
noise, ergodic diffusion processes and Poisson process with Poisson and
self-exciting alternatives. For the first two classes we just show how
Cram\'er-von Mises and Kolmogorov-Smirnov - type tests can be realized using
some known results and for the last models we discuss this problem in
detail. \vspace*{20pt}

\hrule

\section{Diffusion process with small noise}

Suppose that the observed process is the solution of the stochastic
differential equation
\begin{equation}
\label{1}
{\rm d}X_t=S\left(X_t\right)\,{\rm d}t+\varepsilon \,{\rm d}W_t,\qquad
X_0=x_0,\quad 0\leq t\leq T,
\end{equation}
where $W_t,0\leq t\leq T$ is a Wiener process (see, for example, Liptser and
Shiryayev (2001)). We assume that the function $S\left(x\right)$ is two times
continuously differentiable with bounded derivatives. These are not the
minimal conditions for the results presented below, but this assumption
simplifies the exposition.  We are interested in the statistical inference for
this model in the asymptotics of small noise : $\varepsilon \rightarrow
0$. The statistical estimation theory (parametric and nonparametric) was
developed in Kutoyants (1994).

Recall that the stochastic process $X^\varepsilon =\left\{X_t,
0\leq t\leq T\right\}$ converges uniformly in $t\in \left[0,T\right]$ to the
deterministic function $ \left\{x_t, 0\leq t\leq T\right\}$, which is a
solution of the ordinary differential equation
\begin{equation}
\label{2}
\frac{{\rm d}x_t}{{\rm d}t}=S\left(x_t\right),\qquad x_0, \quad 0\leq t\leq T.
\end{equation}
Suppose that the function $S_*\left(x\right)>0 $ for $ x\geq x_0$ and consider
the following problem of hypotheses testing
\begin{eqnarray*}
&&{\scr   H}_0:\qquad \qquad S\left(x\right)=S_*\left(x\right),\quad x_0\leq
x\leq x_T^*\\
&&{\scr H}_1:\qquad \qquad S\left(x\right)\neq S_*\left(x\right),\quad
x_0\leq x\leq x_T^*
\end{eqnarray*}
where we denoted by $x_t^* $ the solution of the equation \eqref{2} under
hypothesis ${\scr H}_0 $:
$$
x_t^*=x_0+\int_{0}^{t}S_*\left(x_v^*\right)\;{\rm d}v,\qquad 0\leq t\leq T.
$$
Hence, we have a simple hypothesis against the composite alternative.

The Cram\'er-von Mises $\left(W_\varepsilon^2 \right)$ and Kolmogorov-Smirnov
$\left(D_\varepsilon \right)$ type statistics for this model of observations
can be
\begin{align*}
W _\varepsilon^2 &=\left[\int_{0}^{T}\frac{{\rm d}t}{S_*\left(x_t^*\right)^2
}\right]^{-2}\;
\int_{0}^{T}\left(\frac{X_t-x_t^*}{\varepsilon\,S_*\left(x_t^*\right)^2
}\right)^2\,{\rm d}t,\\
D_\varepsilon &=\left[\int_{0}^{T}\frac{{\rm d}t}{S_*\left(x_t^*\right)^2
}\right]^{-1/2}\;
\sup_{0\leq t\leq T}\left|\frac{X_t-x_t^*}{S_*\left(x_t^*\right)}\right|.
\end{align*}
It can be shown that these two statistics converge (as $\varepsilon
\rightarrow 0 $) to the following functionals
$$
W_\varepsilon^2 \Longrightarrow \int_{0}^{1}W\left(s\right)^2\,{\rm d}s,
\qquad  \varepsilon ^{-1} D_\varepsilon \Longrightarrow \sup_{0\leq s\leq
1}\left|W\left(s\right)\right|,
$$
where $\left\{W\left(s\right),0\leq s\leq 1\right\}$ is a Wiener process
(see Kutoyants 1994).  Hence the corresponding tests
$$
\psi_\varepsilon\left(X^\varepsilon \right) =1_{\left\{W
_\varepsilon^2>c_\alpha \right\}},\qquad
\phi_\varepsilon\left(X^\varepsilon \right) =1_{\left\{\varepsilon
^{-1}D_\varepsilon >d_\alpha \right\}}
$$
with the constants $c_\alpha , d_\alpha $ defined by the equations
\begin{equation}
\label{3}
\Pb\left\{\int_{0}^{1}W\left(s\right)^2\;{\rm d}s > c_\alpha \right\}=\alpha
,\qquad
\Pb\left\{\sup_{0\leq s\leq 1}\left|W\left(s\right)\right| > d_\alpha
\right\}=\alpha
\end{equation}
are of asymptotic size $\alpha $. Note that the choice of the thresholds
$c_\alpha $ and $d_\alpha $ does not depend on the hypothesis
(distribution-free). This situation is quite close to the classical case
mentioned above.

It is easy to see that if $S\left(x\right)\neq S_*\left(x\right)$, then
$\sup_{0\leq t\leq T}\left|x_t-x_t^*\right|>0$ and $W
_\varepsilon^2\rightarrow \infty $, $ \varepsilon ^{-1}D_\varepsilon
\rightarrow \infty $. Hence these tests are consistent against any fixed
alternative. It is possible to study the power function of this test for local
(contiguous) alternatives of the following form
\begin{equation*}
\label{}
{\rm d}X_t=S_*\left(X_t\right)\,{\rm d}t+\varepsilon \;
\frac{h\left(X_t\right)}{S_*\left(X_t\right)}\;{\rm
d}t+\varepsilon \;{\rm d}W_t,\quad 0\leq t\leq T .
\end{equation*}

We describe the alternatives with the help of the (unknown) function
$h\left(\cdot \right)$. The case $h\left(\cdot \right)\equiv 0$ corresponds to
the hypothesis ${\scr H}_0$. One special class of such nonparametric
alternatives for this model was studied in Iacus and Kutoyants (2001).

Let us introduce the composite (nonparametric)
alternative
$$
{\scr H}_1\quad :\qquad \qquad h\left(\cdot \right)\in
{\cal H}_{\rho },
$$
where
$$
{\cal H}_\rho=\left\{h\left(\cdot \right)\; : \qquad
\int_{x_0}^{x_T}h\left(x\right)^2 \;\mu \left({\rm d}x\right)\geq
\rho \right\} .
$$
To choose alternative we have to precise the ``natural for this problem''
distance described by the measure $\mu \left(\cdot \right)$ and the rate of
$\rho=\rho _\varepsilon $. We show that the choice
$$
\mu \left({\rm d}x\right)=\frac{{\rm d}x}{S_*\left(x\right)^3}
$$
provides for the test statistic the following limit
$$
W _\varepsilon^2 \longrightarrow
\int_{0}^{1}\left[\int_{0}^{s}h_*\left(v\right){\rm d}v+W\left(s\right)
\right]^2 {\rm d}s,
$$
where we denoted
$$
h_*\left(s\right)=u_T^{1/2}h\left(x^*_{u_T s}\right),\qquad
u_T=\int_{0}^{T}\frac{{\rm d}s}{ S_*\left(x_s^*\right)^2}
$$
We see that this problem  is asymptotically equivalent to the   {\sl
signal in white Gaussian noise} problem:
\begin{equation}
\label{swn}
{\rm d}Y_s=h_*\left(s\right)\,{\rm d}s+{\rm d}W\left(s\right),\quad 0\leq
s\leq 1,
\end{equation}
with the Wiener process $W\left(\cdot \right)$. It is easy to see that even
for fixed $\rho >0$ without further restrictions on the smoothness of the
function $h_*\left(\cdot \right)$ the {\sl uniformly good} testing is
impossible. For example, if we put
$$
h_n\left(x\right)=c\,S_*\left(x\right)^3\;\cos\left[n\left(x-x_0\right)\right]
$$
then for the power function of the test we have
$$
\inf_{h\left(\cdot \right)\in {\cal H}_\rho} \beta \left(\psi_\varepsilon
,h\right)\leq \beta \left(\psi_\varepsilon ,h_n\right) \longrightarrow \alpha .
$$
The details can be found in Kutoyants (2006). The construction of the
uniformly consistent tests requires a different approach (see Ingster and
Suslina (2003)).

Note as well that if the diffusion process is
$$
{\rm d}X_t=S\left(X_t\right)\,{\rm d}t+ \varepsilon \sigma
\left(X_t\right)\;{\rm d}W_t, \qquad X_0=x_0,\qquad 0\leq t\leq T ,
$$
then we can put
$$
W _\varepsilon^2 =\left[\int_{0}^{T}\left(\frac{\sigma
\left(x_t^*\right)}{S_*\left(x_t^*\right)
}\; \right)^2  {\rm d}t\right]^{-2}\;
\int_{0}^{T}\left(\frac{X_t-x_t^*}{\varepsilon\,S_*\left(x_t^*\right)^2
}\right)^2\,{\rm d}t
$$
and have the same results as above (see Kutoyants (2006)).

\vspace*{20pt}

\hrule
\section{Ergodic diffusion processes}

Suppose that the observed process is one dimensional diffusion process
\begin{equation}
\label{4}
{\rm d}X_t=S\left(X_t\right)\,{\rm d}t+ {\rm d}W_t,\qquad X_0,\qquad 0\leq
t\leq T ,
\end{equation}
where the trend coefficient $S\left(x\right)$ satisfies the conditions of the
existence and uniqueness of the solution of this equation and this solution has
ergodic properties, i.e., there exists an invariant probability distribution
$F_S\left(x\right)$, and for any integrable w.r.t. this distribution function
$g\left(x\right)$ the law of large numbers  holds
$$
\frac{1}{T}\int_{0}^{T}g\left(X_t\right)\;{\rm d}t\longrightarrow
\int_{-\infty }^{\infty }g\left(x\right) \; {\rm d}F_S\left(x\right).
$$
These conditions can be found, for example, in Kutoyants (2004).

Recall that the invariant   density function
$f_S\left(x\right)$  is defined  by the equality
$$
f_S\left(x\right)=G\left(S\right)^{-1}\exp\left\{2\int_{0}^{x}S
\left(y\right)\;{\rm d}y\right\},
$$
where $G\left(S\right)$ is the normalising constant.

We consider two types of tests. The first one is a direct analogue of the
classical Cram\'er-von Mises and Kolmogorov-Smirnov tests based on empirical
distribution and density functions and the second follows the considered above
(small noise)  construction of tests.

The invariant distribution function $F_S\left(x\right)$ and this density
function can be estimated by the {\sl empirical distribution function} $\hat
F_T\left(x\right) $ and by the {\sl local time type} estimator $\hat
f_T\left(x\right)$ defined by the equalities
$$
\hat F_T\left(x\right)=\frac{1}{T}\int_{0}^{T}1_{\left\{X_t<x\right\}}\;{\rm
d}t,\qquad \quad \hat
f_T\left(x\right)=\frac{2}{T}\int_{0}^{T}1_{\left\{X_t<x\right\}}\;{\rm d}X_t
$$
respectively. Note that both of them are unbiased:
$$
\Ex_S\hat F_T\left(x\right)=
F_S\left(x\right) ,\qquad \Ex_S\hat f_T\left(x\right)=
f_S\left(x\right),
$$
admit the representations
\begin{align*}
\eta _T\left(x\right)&=-\frac{2}{\sqrt{T}}\int_{0}^{T}
\frac{F_S\left(X_t\wedge x
\right)-F_S\left(X_t\right)F_S\left(x\right)}{f_S\left(X_t\right)}\;{\rm
d}W_t+o\left(1\right),\\ \zeta
_T\left(x\right)&=-\frac{2f_S\left(x\right)}{\sqrt{T}}\int_{0}^{T}
\frac{1_{\left\{X_t
>x\right\}}-F_S\left(X_t \right) }{f_S\left(X_t \right) }\;{\rm
d}W_t+o\left(1\right)
\end{align*}
and are $\sqrt{T}$ asymptotically normal (as $T\rightarrow \infty $)
\begin{align*}
\eta _T\left(x\right)=\sqrt{T}\left(\hat F_T\left(x\right)-F_S\left(x\right)
\right)& \Longrightarrow {\cal N}\left(0, d_F\left(S,x\right)^2 \right),\\
\zeta _T\left(x\right)=\sqrt{T}\left(\hat f_T\left(x\right)-f_S\left(x\right)
\right)&\Longrightarrow {\cal N}\left(0, d_f\left(S,x\right)^2 \right).
\end{align*}

Let us fix a simple  (basic) hypothesis
$$
{\scr H}_0\quad :\qquad S\left(x\right)\equiv S_*\left(x\right).
$$
Then to test this hypothesis we can use these estimators for construction of
the Cram\'er-von Mises and Kolmogorov-Smirnov type test statistics
\begin{align*}
W_T^2&=T\int_{-\infty }^{\infty }\left[\hat
F_T\left(x\right)-F_{S_*}\left(x\right) \right]^2\; {\rm
d}F_{S_*}\left(x\right),\\
D_T&=\sup_x\left|\hat F_T\left(x\right)-F_{S_*}\left(x\right)\right|
\end{align*}
and
\begin{align*}
V_T^2&=T\int_{-\infty }^{\infty }\left[\hat
f_T\left(x\right)-f_{S_*}\left(x\right) \right]^2\; {\rm
d}F_{S_*}\left(x\right),\\
d_T&=\sup_x\left|\hat f_T\left(x\right)-f_{S_*}\left(x\right)\right|
\end{align*}
respectively. Unfortunately, all these statistics are not distribution-free
even asymptotically and the choice of the corresponding thresholds for the
tests is much more complicated. Indeed, it was shown that the random functions
$\left(\eta _T\left(x\right), x\in R \right)$ and $\left(\zeta
_T\left(x\right), x\in R \right)$ converge in the space $\left({\scr C}_0,
{\goth B}\right)$ (of continuous functions decreasing to zero at infinity) to
the zero mean Gaussian processes $\left(\eta \left(x\right), x\in R \right)$
and $\left(\zeta \left(x\right), x\in R \right)$ respectively with the
covariance functions (we omit the index $S_*$ of functions
$f_{S_*}\left(x\right) $ and $F_{S_*}\left(x\right) $ below)
\begin{align*}
R_F\left(x,y\right)&=\Ex_{S_*}\left[\eta \left(x\right)\eta
\left(y\right)\right]\\
&=4\Ex_{S_*}\left(\frac{\left[F\left(\xi\wedge x
\right)-F\left(\xi\right)F\left(x\right)\right]
\left[F\left(\xi\wedge y
\right)-F\left(\xi \right)F\left(y\right)\right] }{f\left(\xi \right)^2
}\right)\\
R_f\left(x,y\right)&=\Ex_{S_*}\left[\zeta  \left(x\right)\zeta
\left(y\right)\right]\\
&=
4f\left(x\right)f\left(y\right) \Ex_{S_*}\left(\frac{\left[1_{\left\{\xi
>x\right\}}-F\left(\xi \right)\right]\left[1_{\left\{\xi
>y\right\}}-F\left(\xi  \right)\right] }{f\left(\xi  \right)^2 }\right).
\end{align*}
Here $\xi $ is a random variable with the distribution function
$F_{S_*}\left(x\right)
$. Of course,
$$
d_F\left(S,x\right)^2=\Ex_S\left[\eta \left(x\right)^2\right],\qquad \qquad
d_f\left(S,x\right)^2=\Ex_S\left[\zeta  \left(x\right)^2\right].
$$
Using this weak convergence it is shown that these statistics converge in
distribution (under hypothesis) to the following limits (as $T\rightarrow
\infty $)
\begin{align*}
&W_T^2\Longrightarrow \int_{-\infty }^{\infty }\eta \left(x\right)^2\;{\rm
d}F_{S_*}\left(x\right),\qquad
T^{1/2}D_T\Longrightarrow \sup_x\left|\eta \left(x\right)\right|,\\
&V_T^2\Longrightarrow \int_{-\infty }^{\infty }\zeta  \left(x\right)^2\;{\rm
d}F_{S_*}\left(x\right),\qquad
T^{1/2}d_T\Longrightarrow \sup_x\left|\zeta \left(x\right)\right|.
\end{align*}
The conditions and the proofs of all these properties can be found in
Kutoyants (2004), where essentially different statistical problems were
studied, but the calculus are quite close to what we need here.

Note that the Kolmogorov-Smirnov test for ergodic diffusion was studied in
Fournie (1992) (see as well Fournie and Kutoyants (1993) for further details),
and the weak convergence of the process $\eta _T\left(\cdot \right)$ was
obtained in Negri (1998).

The Cram\'er-von Mises and Kolmogorov-Smirnov type tests based on these
statistics  are
\begin{align*}
\Psi_T\left(X^T\right)&=1_{\left\{W_T^2>C_\alpha \right\}},\qquad
\Phi_T\left(X^T\right)=1_{\left\{T^{1/2}D_T>D_\alpha  \right\}},\\
\psi_T\left(X^T\right)&=1_{\left\{V_T^2>c_\alpha \right\}},\qquad
\phi_T\left(X^T\right)=1_{\left\{T^{1/2}d_T>d_\alpha  \right\}}
\end{align*}
with appropriate constants.

The contiguous alternatives can be introduced by the following way
$$
S\left(x\right)=S_*\left(x\right)+\frac{h\left(x\right)}{\sqrt{T}}.
$$
Then we obtain for the Cram\'er-von Mises  statistics the limits (see, Kutoyants
(2004))
\begin{align*}
&W_T^2\Longrightarrow \int_{-\infty }^{\infty }\left[2\Ex_{S_*}
\left(\left[1_{\left\{\xi <x\right\}}-F_{S_*}\left(x\varphi \right)\right]
\int_{0
}^{\xi}h\left(s\right)\;{\rm d}s\right)+\eta \left(x\right)\right]^2\;{\rm
d}F_{S_*}\left(x\right),\\
&V_T^2\Longrightarrow \int_{-\infty }^{\infty }\left[2f_{S_*}\left(x\right)\, \Ex_{S_*}\int_{\xi
}^{x}h\left(s\right)\;{\rm d}s+\zeta  \left(x\right)\right]^2\;{\rm
d}F_{S_*}\left(x\right).
\end{align*}

Note that the transformation $Y_t=F_{S_*}\left(X_t\right)$ simplifies the
writing, because the diffusion process $Y_t$ satisfies the differential
equation
$$
{\rm
d}Y_t=f_{S_*}\left(X_t\right)\left[2S_*\left(X_t\right){\rm
d}t+ \;{\rm d}W_t\right],\qquad Y_0=F_{S_*}\left(X_0\right)
$$
with reflecting bounds in 0 and 1 and (under hypothesis) has uniform on
$\left[0,1\right]$ invariant distribution. Therefore,
$$
W_T^2\Longrightarrow \int_{0}^{1}V\left(s\right)^2{\rm d}s,\qquad
T^{1/2}D_T\Longrightarrow \sup_{0\leq s\leq 1} \left|V\left(s\right)\right|,
$$
but the covariance structure of the Gaussian process $\left\{V\left(s\right),
0\leq s\leq 1 \right\}$ can be quite complicated.

\bigskip

To obtain asymptotically  distribution-free Cram\'er-von Mises type test we can
use another statistic, which is similar to that of the preceding
section. Let us introduce
$$
\tilde W_T^2=\frac{1}{T^2}\int_{0}^{T}\left[X_t-X_0-\int_{0}^{t}
S_*\left(X_v\right)\,{\rm d}v\right]^2{\rm d}t.
$$
Then we have immediately (under hypothesis)
$$
\tilde W_T^2=\frac{1}{T^2}\int_{0}^{T}W_t^2\;{\rm
d}t=\int_{0}^{1}W\left(s\right)^2\;{\rm d}s ,
$$
where we put $t=sT$ and $W\left(s\right)=T^{-1/2}W_{sT}$. Under alternative we
have
\begin{align*}
\tilde W_T^2&=\frac{1}{T^2}\int_{0}^{T}\left[W_t+\frac{1}{\sqrt{T}}
\int_{0}^{t}h\left(X_v \right){\rm d}v\right]^2{\rm d}t\\
&=
\frac{1}{T}\int_{0}^{T}\left[\frac{W_t}{\sqrt{T}}+\frac{t}{T} \,\frac{1}{t}
\int_{0}^{t}h\left(X_v \right){\rm d}v\right]^2{\rm d}t.
\end{align*}
The stochastic process $X_t$ is ergodic, hence
$$
\frac{1}{t}\int_{0}^{t}h\left(X_v \right){\rm d}v\longrightarrow
\Ex_{S_*}h\left(\xi\right)=\int_{-\infty }^{\infty
}h\left(x\right)f_{S_*}\left(x\right){\rm d}x \equiv \rho_h
$$
as $t\rightarrow \infty $. It can be shown (see section 2.3 in Kutoyants
(2004), where we have the similar calculus in another problem)
that
$$
\tilde W_T^2\Longrightarrow \int_{0}^{1}\left[\rho_h
\,s+W\left(s\right)\right]^2{\rm d}s.
$$

Therefore the power function of the test
$\psi\left(X^T\right)=1_{\left\{\tilde W_T^2>c_\alpha \right\}}$ converges to
the function
$$
\beta_\psi\left(\rho_h \right)=\Pb\left( \int_{0}^{1}\left[\rho_h
\,s+W\left(s\right)\right]^2{\rm d}s>c_\alpha \right) .
$$

Using standard calculus  we can show that for the corresponding
Kol\-mo\-go\-rov-Smirnov type test the limit will be
$$
\beta_\phi\left(\rho_h \right)=\Pb\left( \sup_{0\leq s\leq 1}\left|\rho_h
\,s+W\left(s\right)\right|>c_\alpha \right) .
$$

These two limit power functions are the same as in the next section devoted to
self-exciting alternatives of the Poisson process. We calculate these
functions with the help of simulations in Section 5 below.

Note that if the diffusion process is
$$
{\rm d}X_t=S\left(X_t\right)\,{\rm d}t+  \sigma \left(X_t\right)\;{\rm d}W_t,
\qquad X_0,\qquad 0\leq t\leq T ,
$$
but the functions $S\left(\cdot \right)$ and  $\sigma \left(\cdot \right)$ are
such that the process is  ergodic then we introduce the statistics
$$
\hat W_T^2=\frac{1}{T^2\,\Ex_{S_*}\left[\sigma \left(\xi \right)^2\right]}\;
\int_{0}^{T}\left[X_t-X_0-\int_{0}^{t} S_*\left(X_v\right)\;{\rm d}v
\right]^2{\rm d}t .
$$
Here $\xi $ is random variable with the invariant density function
$$
f_{S_*}\left(x\right)=\frac{1}{G\left(S_*\right)\sigma
\left(x\right)^2}\;\exp\left\{2\int_{0}^{x }\frac{S_*\left(y\right)}{\sigma
\left(y\right)^2}\;{\rm d}y\right\} .
$$
This statistic under hypothesis is equal to
\begin{align*}
\hat W_T^2&=\frac{1}{T^2\,\Ex_{S_*}\left[\sigma \left(\xi \right)^2\right]}\;
\int_{0}^{T}\left[\int_{0}^{t} \sigma \left(X_v\right){\rm d}W_v
\right]^2{\rm d}t  \\
&=\frac{1}{T\,\Ex_{S_*}\left[\sigma \left(\xi
\right)^2\right]}\;\int_{0}^{T}\left[
\frac{1}{\sqrt{T}} \int_{0}^{t} \sigma \left(X_v\right){\rm d}W_v
\right]^2{\rm d}t  .
\end{align*}
The   stochastic integral by the central limit theorem is asymptotically
normal
$$
\eta _t=\frac{1}{\sqrt{t \Ex_{S_*}\left[\sigma \left(\xi \right)^2\right] }
}\int_{0}^{t} \sigma \left(X_v\right){\rm d}W_v\Longrightarrow {\cal
N}\left(0,1 \right)
$$
and moreover it can be shown that the vector of such integrals converges in
distribution to the Wiener process
$$
\Bigl(\eta _{s_1T},\ldots,\eta _{s_kT}\Bigr)\Longrightarrow
\left(W\left(s_1\right),\ldots,W\left(s_k\right) \right)
$$
for any finite collection of $0\leq s_1<s_2<\ldots <s_k\leq 1$. Therefore,
under mild regularity conditions it can be proved that
$$
\hat W_T^2\Longrightarrow \int_{0}^{1}W\left(s\right)^2\;{\rm d}s.
$$
The power function has the same limit,
$$
\beta_\psi\left(\rho_h \right)=\Pb\left( \int_{0}^{1}\left[\rho_h
\,s+W\left(s\right)\right]^2{\rm d}s>c_\alpha \right) .
$$
 but with
$$
\rho _h=\frac{\Ex_{S_*}h\left(\xi \right)}{\sqrt{\Ex_{S_*} \left[\sigma
\left(\xi \right)^2\right]  }}.
$$
The similar consideration can be done for the Kolmogorov-Smirnov type test
too.

We see that both tests can not distinguish the alternatives with $h\left(\cdot
\right)$ such that $\Ex_{S_*}h\left(\xi \right)=0 $. Note that for ergodic
processes usually we have $\Ex_{S}S\left(\xi \right)=0$ and
$\Ex_{S_*+h/\sqrt{T}}\left[S_*\left(\xi \right) +T^{-1/2}h\left(\xi
\right)\right]=0$ with corresponding random variables $\xi $, but this does
not imply $\Ex_{S_*}h\left(\xi \right)=0 $.

\vspace*{20pt}

\hrule
\section{Poisson and self-exciting processes}

Poisson process is one of the simplest point processes and before taking any
other model it is useful first of all to check the hypothesis the observed
sequence of events, say, $0<t_1,\ldots,t_N<T$ corresponds to a Poisson
process. It is natural in many problems to suppose that this Poisson process
is periodic of known period. For example, many daily events, signal
transmission in optical communication, season variations etc. Another model of
point processes as well frequently used is self-exciting stationary point
process introduced in Hawkes (1972).  As any stationary process it can as well
describe the periodic changes due to the particular form of its spectral
density.

Recall that for the Poisson process $X_t,t\geq 0$ of intensity function
$S\left(t\right), t\geq 0 $ we have ($X_t$ is the counting process)
$$
\Pb\left\{X_t-X_s=k\right\}=\left(k!\right)^{-1}\,\left(\Lambda
\left(t\right)-\Lambda \left(s\right)\right)^k \;\exp\left\{\Lambda
\left(s\right)-\Lambda \left(t\right)\right\} ,
$$
where we suppose that $s<t$ and put
$$
\Lambda \left(t\right)=\int_{0}^{t}S\left(v\right)\,{\rm d}v.
$$
The self-exciting process $X_t,t\geq 0$  admits the representation
$$
X_t=\int_{0}^{t}S\left(s,X\right)\;{\rm d}s+ \pi _t,
$$
where $\pi _t,t\geq 0$ is local martingale and the  intensity function
$$
S\left(t,X\right)=S+\int_{0}^{t}g\left(t-s\right)\;{\rm
d}X_s=S+\sum_{t_i<T}^{}g\left(t-t_i\right).
$$
It is supposed that
$$
\rho =\int_{0}^{\infty }g\left(t\right)\;{\rm d}t<1.
$$
Under this condition the self-exciting process is a stationary point process
with the rate
$$
\mu =\frac{S}{1-\rho }
$$
and the spectral density
$$
f\left(\lambda \right)=\frac{\mu }{2\pi
\left|1-G\left(\lambda \right)\right|^2},\qquad G\left(\lambda \right)
=\int_{0}^{\infty }e^{{\rm i}\lambda t}g\left(t\right)\;{\rm d}t
$$
(see Hawkes (1972) or Daley and Vere-Jones (2003) for details).

We consider two problems: Poisson against another Poisson and Poisson against
a close self-exciting point process.
The first one is to test the simple (basic)  hypothesis
$$
{\scr H}_0\qquad :\qquad S\left(t\right)\equiv S_*\left(t\right),\quad t\geq 0
$$
where $S_*\left(t\right)$ is known periodic function of period $\tau $,
against the composite alternative
$$
{\scr H}_1\qquad :\qquad S\left(t\right)\neq S_*\left(t\right),\quad t\geq 0,
$$
but $S\left(t\right)$ is always $\tau $-periodic.

Let us denote $X_j\left(t\right)=X_{\tau \left(j-1\right)+t}-X_{\tau
\left(j-1\right)}$, $j=1,\ldots,n$, suppose that $T=n\tau $ and put
$$
\hat\Lambda _n\left(t\right)=\frac{1}{n}\sum_{j=1}^{n}X_j\left(t\right).
$$

The corresponding goodness-of-fit tests of Cram\'er-von Mises and
Kol\-mo\-go\-rov-Smirnov type can be based on the statistics
\begin{align*}
&W_n^2=\Lambda_* \left(\tau \right)^{-2}n\int_{0}^{\tau }\left[\hat\Lambda
_n\left(t\right)-\Lambda_* \left(t\right)\right]^2{\rm d}\Lambda_*
\left(t\right),\\ & D_n=\Lambda_* \left(\tau \right)^{-1/2}\sup_{0\leq t\leq
\tau } \left|\hat\Lambda _n\left(t\right)-\Lambda_* \left(t\right)\right|.
\end{align*}
It can be shown that
$$
W_n^2\Longrightarrow \int_{0}^{1}W\left(s\right)^2{\rm d}s,\qquad \sqrt{n}\;
D_n\Longrightarrow \sup_{0\leq s\leq 1 }\left|W\left(s\right)\right|
$$
where $\left\{W\left(s\right),0\leq s\leq 1\right\}$ is a Wiener process (see
Kutoyants (1998)). Hence these statistics are asymptotically distribution-free
and the tests
$$
\psi_n\left(X^T\right)=1_{\left\{W_n^2>c_\alpha \right\}},\qquad \quad
\phi_n\left(X^T\right)=1_{\left\{\sqrt{n}D_n>d_\alpha \right\}}
$$
with the constants $c_\alpha , d_\alpha $ taken from the equations \eqref{3},
are of asymptotic size $\alpha $.

Let us describe the close contiguous alternatives which reduce asymptotically
this problem to {\sl signal in white Gaussian noise} model \eqref{swn}. We put
$$
\Lambda \left(t\right)=\Lambda _*\left(t\right)+\frac{1}{\sqrt{n\Lambda
_*\left(\tau \right) }}\int_{0}^{t}h\left(u\left(v\right) \right) \;{\rm d}
\Lambda _*\left(v\right),\qquad u\left(v\right)=\frac{\Lambda _*\left(v\right)
}{\Lambda _*\left(\tau \right)}.
$$
Here $h\left(\cdot \right)$ is an arbitrary function defining the alternative.
Then if $\Lambda \left(t\right)$ satisfies this equality we have the
convergence
$$
W_n^2\Longrightarrow \int_{0}^{1}\left[\int_{0}^{s}h\left(v\right){\rm
d}v+W\left(s\right)\right]^2{\rm d}s.
$$

This convergence describes the power function of the Cram\'er-von Mises type
test under these alternatives.

\bigskip

The second problem is to test the hypothesis
$$
{\scr H}_0\qquad :\qquad S\left(t\right)= S_*,\quad t\geq 0
$$
against nonparametric close (contiguous) alternative
$$
{\scr H}_1\qquad :\qquad S\left(t\right)=
S_*+\frac{1}{\sqrt{T}}\int_{0}^{t}h\left(t-s\right){\rm d}X_t,\quad t\geq 0 ,
$$
We consider the alternatives with the functions $h\left(\cdot \right)\geq 0$
having compact support and bounded.

We have $\Lambda _*\left(t\right)=S_*\;t$ and for some fixed $\tau >0$ we can
construct the same statistics
$$
W_n^2=\frac{n}{S_*\tau^2}\int_{0}^{\tau }\left[\hat\Lambda
_n\left(t\right)-S_*\;t\right]^2\;{\rm d}t,\quad
D_n=\left(S_*\;\tau \right)^{-1/2}\sup_{0\leq t\leq \tau }\left|\hat\Lambda
_n\left(t\right)-S_*\;t \right|.
$$
Of course, they have the same limits under hypothesis
$$
W_n^2\Longrightarrow \int_{0}^{1}W\left(s\right)^2{\rm d}s,\qquad
\sqrt{n}D_n\Longrightarrow \sup_{0\leq s\leq 1 } \left|W\left(s\right)\right|.
$$

To describe their behaviour under any
fixed alternative $h\left(\cdot \right)$ we have to find the limit
distribution of the vector
$$
{\bf w}_n=\bigl(w_n\left(t_1\bigr),\ldots,w_n\left(t_k\right)\right),\qquad
w_n\left(t_l\right)= \frac{1}{\sqrt{S_* \tau\;
n}}\sum_{j=1}^{n}\left[X_j\left(t_l\right)-S_*t_l\right],
$$
where $0\leq t_l \leq \tau $. We know that this vector under
hypothesis is asymptotically normal
$$
{\cal L}_0\left\{{\bf w}_n \right\}\Longrightarrow {\cal N}\left({\bf 0}, {\bf
R}\right)
$$
with covariance matrix
$$
{\bf R}=\left(R_{lm}\right)_{k\times k}, \qquad
R_{lm}=\tau ^{-1}\min\left(t_l,t_m\right).
$$
Moreover, it was shown in Dachian and Kutoyants (2006) that for such
alternatives the likelihood ratio is locally asymptotically normal, i.e., the
likelihood ratio admits the representation
$$
Z_n\left(h\right)=\exp\left\{\Delta _n\left(h,X^n\right)-\frac{1}{2}\,{\rm
I}\left(h \right)+r_n\left(h,X^n\right)\right\}
$$
where
\begin{align*}
&\Delta _n\left(h,X^n\right)=\frac{1}{S_*\sqrt{\tau n}}\int_{0}^{\tau
n}\int_{0}^{t-}h\left(t-s\right) \,{\rm d}X_s\;\left[{\rm d}X_t-S_*{\rm
d}t\right],\\
&{\rm I}\left(h \right)=\int_{0}^{\infty }h\left(t\right)^2{\rm
d}t+S_*\left(\int_{0}^{\infty }h\left(t\right){\rm d}t\right)^2
\end{align*}
and
\begin{equation}
\label{7}
\Delta _n\left(h,X^n\right)\Longrightarrow    {\cal N}\left(0,{\rm I}\left(h
\right) \right),\qquad         r_n\left(h,X^n\right)\rightarrow 0.
\end{equation}
To use the Third Le Cam's Lemma we describe the limit behaviour of the vector
$\left(\Delta _n\left(h,X^n\right),{\bf w}_n \right)$. For the covariance
${\bf Q}=\left( Q_{lm}\right), l,m=0,1,\ldots,k$ of this vector we have
$$
\Ex_0\Delta _n\left(h,X^n\right)=0,\qquad  { Q}_{00}=\Ex_0\Delta
_n\left(h,X^n\right)^2={\rm I}\left(h \right)\left(1+o\left(1 \right)\right).
$$
Further, let us denote ${\rm d}\pi _t={\rm d}X_t-S_*{\rm d}t $ and
$H\left(t\right)=\int_{0}^{t-}h\left(t-s\right) \,{\rm d}X_s $, then we can
write
\begin{align*}
&{ Q}_{0l}=\Ex_0\left[\Delta _n\left(h,X^n\right)w_n\left(t_l\right)\right]\\
&\quad = \frac{1}{n
S_*^{3/2}\tau }\;\Ex_0\left( \sum_{j=1}^{n}\int_{\tau
\left(j-1\right)}^{\tau j}H\left(t\right){\rm d}\pi
_t \,  \,\sum_{i=1}^{n}\int_{\tau \left(i-1\right)}^{\tau\left(i-1\right)+t_l}
{\rm d}\pi _t\right)\\
&\quad =\frac{1}{n\tau \sqrt{S_* }}\sum_{j=1}^{n}\int_{\tau
\left(j-1\right)}^{\tau\left(j-1\right)+t_l}\Ex_0 H\left(t\right)\; {\rm d}t=
\frac{t_l}{\tau }\;\sqrt{S_*}\;\int_{0}^{\infty }h\left(t\right){\rm
d}t\;\left(1+o\left(1 \right)\right),
\end{align*}
because
$$
\Ex_0 H\left(t\right)=S_*\;\int_{0}^{t-}h\left(t-s\right){\rm
d}s= S_*\;\int_{0}^{\infty
}h\left(s\right){\rm d}s
$$
for the large values of $t$ (such that $\left[0,t\right]$ covers the support
of $h\left(\cdot \right)$).

Therefore, if we denote
$$
\bar h=\int_{0}^{\infty }h\left(s\right){\rm d}s
$$
then
$$
{ Q}_{0l}={ Q}_{l0}=\frac{t_l}{\tau }\;\sqrt{S_*}\;\bar h.
$$

The proof of the Theorem 1 in Dachian and Kutoyants (2006) can be applied to
the linear combination of $ \Delta _n\left(h,X^n\right)$ and
$w_n\left(t_1\right),\ldots,w_n\left(t_k\right)$ and this yields the
asymptotic normality
$$
{\cal L}_0\Bigl(\Delta _n\left(h,X^n\right),{\bf w}_n \Bigr)\Longrightarrow
{\cal N}\left({\bf 0},{\bf Q}\right).
$$
Hence by the Third Lemma of Le Cam we obtain the asymptotic normality of the
vector ${\bf w}_n$
$$
{\cal L}_h\Bigl({\bf w}_n \Bigr)\Longrightarrow {\cal
L}\left(W\left(s_1\right)+s_1\,\sqrt{S_*}\;\bar
h,\ldots,W\left(s_k\right)+s_k\,\sqrt{S_*}\;\bar h\right),
$$
where we put $t_l=\tau \,s_l$. This weak convergence together with the
estimates like
$$
\Ex_h\left|w_n\left(t_1\right)-w_n\left(t_2\right)\right|^2\leq
C\left|t_1-t_2\right|
$$
provides the convergence (under alternative)
$$
W_n^2\Longrightarrow \int_{0}^{1}\left[\sqrt{S_*}\;\bar
h\,s+W\left(s\right)\right]^2{\rm d}s.
$$

We see that the limit experiment is of the type
$$
dY_s=\sqrt{S_*}\;\bar h\,{\rm d}s+{\rm d}W\left(s\right),\quad Y_0=0,\quad
0\leq s\leq 1.
$$

The power $\beta(\psi_n,h)$ of the Cramer-von Mises type test
$\psi_n(X^n)=1_{\{W_n^2>c_\alpha\}}$ is a function of the real parameter
$\rho_h =\sqrt{S_*}\;\bar h$
$$
\beta \left(W_n,h\right)=
\Pb\left(\int_{0}^{1}\left[\rho_h \,s+W\left(s\right)\right]^2{\rm d}s>c_\alpha
\right)+o\left(1 \right)=\beta_\psi \left(\rho_h \right)+o\left(1 \right).
$$

Using the arguments of Lemma 6.2 in Kutoyants (1998) it can be shown that for
the Kolmogorov-Smirnov type test we have the convergence
$$
\sqrt{n}D_n\Longrightarrow \sup_{0\leq s\leq
1}\left|\rho_h\,s+W\left(s\right)\right|.
$$
The limit power function is
$$
\beta_\phi \left(\rho_h \right)=\Pb\left(\sup_{0\leq s\leq 1}\left|\rho_h
\,s+W\left(s\right)\right|>d_\alpha\right) .
$$

These two limit power functions will be obtained by simulation in the next
section.

\vspace*{20pt}

\hrule

\section{Simulation}

First, we present the simulation of the thresholds $c_\alpha$ and $d_\alpha$
of our Cra\-m\'er-von Mises and Kolmogorov-Smirnov type tests. Since these
thresholds are given by the equations~\eqref{3}, we obtain them by simulating
$10^7$ trajectories of a Wiener process on $[0{,}1]$ and calculating empirical
$1-\alpha$ quantiles of the statistics
$$
W^2=\int_0^1W(s)^2\;{\rm d}s\qquad\text{and}\qquad D=
\sup_{0\leq s\leq 1}\left|W\left(s\right)\right|
$$
respectively. Note that the distribution of $W^2$ coincides with the
distribution of the quadratic form
$$
W^2=\sum_{k=1}^{\infty }\frac{\zeta _k^2}{\left(\pi k\right)^2} ,\qquad \zeta
_k \: \: {\rm i.i.d.}\; \; \sim {\cal N}\left(0,1\right)
$$
and both distributions are extensively studied (see (1.9.4(1)) and (1.15.4) in
Borodin and Salmienen (2002)). The analytical expressions are quite
complicated and we would like to compare by simulation $c_\alpha$ and
$d_\alpha$ with the real (finite time) thresholds giving the tests of exact
size $\alpha$, that is $c_\alpha^T$ and $d_\alpha^T$ given by equations
$$
\Pb\left\{W_n^2 > c_\alpha^T \right\}=\alpha
\qquad\text{and}\qquad
\Pb\left\{\sqrt{n}D_n > d_\alpha^T \right\}=\alpha
$$
respectively. We choose $S^*=1$ and obtain $c_\alpha^T$ and $d_\alpha^T$ by
simulating $10^7$ trajectories of a Poisson process of intensity $1$ on
$[0{,}T]$ and calculating empirical $1-\alpha$ quantiles of the statistics
$W_n^2$ and $\sqrt{n}D_n$. The thresholds simulated for $T=10$, $T=100$ and
for the limiting case are presented in Fig.~1. The lower curves correspond to
the Cram\'er-von Mises type test, and the upper ones to the Kolmogorov-Smirnov
type test. As we can see, for $T=100$ the real thresholds are already
indistinguishable from the limiting ones, especially in the case of the
Cram\'er-von Mises type test.

\onepic{TshNB.eps}{Fig.~1: Threshold choice}

It is interesting to compare the asymptotics of the Cram\'er-von Mises and
Kolmogorov-Smir\-nov type tests with the locally asymptotically uniformly most
powerful (LAUMP) test
$$
\hat\phi_n\left(X^n\right)=1_{\left\{\delta _T>z_\alpha \right\}},\qquad
\delta _T=\frac{X_{n\tau }-S_*n\tau }{\sqrt{S_*n\tau }}
$$
proposed for this problem in Dachian and Kutoyants (2006). Here $z_\alpha $ is
$1-\alpha $ quantile of the standard Gaussian law, $\Pb\left(\zeta >z_\alpha
\right)=\alpha $, $\zeta \sim {\cal N}\left(0,1\right)$. The limit power
function of $\hat\phi_n$ is
$$
\beta_{\hat\phi} \left(\rho_h\right)=\Pb\left(\rho_h+\zeta >z_\alpha
\right).
$$
In Fig.~2 we compare the limit power functions $\beta_{\psi} \left(\rho
\right),\beta_{\phi} \left(\rho \right) $ and $\beta_{\hat\phi} \left(\rho
\right)$. The last one can clearly be calculated directly, and the first two
are obtained by simulating $10^7$ trajectories of a Wiener process on
$[0{,}1]$ and calculating empirical frequencies of the events
$$
\left\{\int_{0}^{1}\left[\rho\,s+W\left(s\right)\right]^2{\rm d}s>c_\alpha
\right\}
\qquad\text{and}\qquad
\left\{\sup_{0\leq s\leq 1}\left|\rho\,s+W\left(s\right)\right|>d_\alpha
\right\}
$$
respectively.

\onepic{LimPowNB.eps}{Fig.~2: Limit power functions}

The simulation shows the exact (quantitative) comparison of the limit power
functions. We see that the power of LAUMP test is higher that the two others
and this is of course evident. We see also that the Kolmogorov-Smirnov type
test is more powerful that the Cram\'er-von Mises type test.

\vspace*{20pt}

\hrule

\end{document}